# FINITE ELEMENT-BASED STRUCTURAL OPTIMIZATION OF LARGE SYSTEM MODELS UNDER BUCKLING CONSTRAINTS

**Alexander Lavin**[1]
Carnegie Mellon University
Pittsburgh, PA
alavin@alumni.cmu.edu

**Kenjji Shimada**[2]
Carnegie Mellon University
Pittsburgh, PA
shimada@cmu.edu

**Giovanni Greco**[3]
Blue Origin, LLC
Seattle WA
ggreco@blueorigin.com

---

[1] Master's student, Dept. of Mechanical Engineering, Carnegie Institute of Technology, 2014
[2] Professor of Mechanical Engineering, Director of CERLAB, Carnegie Institute of Technology
[3] Product Group Lead, Structural Configuration and Design, Blue Origin, LLC

## ABSTRACT

Optimization of large structures of multiple components is essential to many industries for minimizing mass, especially the design of aerospace vehicles. Optimizing a single primary load member independently of all other primary structures is an incomplete process, due to the redistribution of internal loads, as the stiffness distribution changes. That is, optimizing a component changes joint loads, which then calls for a new optimization – changing internal loads changes the optimum. This is particularly evident under buckling (stability) constraints. The goal is to develop a finite element-based optimization approach which can be used to optimize each component of a large, primary structure assembly. The optimization objective function will be to minimize mass for the system, and the constraints will be both stress constraints as well as buckling constraints. The research aims to improve both the solution and practical usability of these models. The system of interest is a spacecraft fuselage, of which the member components are panels throughout the structure. We present analyses of several main optimization methods, and define a new algorithm to solve this problem, *eigenOpt*.

## I. INTRODUCTION

In the aerospace industry, where fuel budgets are a priority, optimization is used in vehicle design to minimize mass. The use of finite element analysis (FEA) -based optimization has grown, and with it the size of problems that optimization engines are being called to address. With large numbers of design variables and constraints, typical methods are inefficient. Gradient-based methods, the standard for optimization problems, are theoretically capable of handling large problems. But a main issue arises with the core feature of these algorithms: solving a computationally expensive sub-optimization task to find the search direction. Optimization over a large design space becomes computationally intense, preventing convergence on a solution. Non gradient-based methods, such as genetic algorithms, have proven unreliable and inefficient for problems involving more than a few variables [1]. Across optimization methods, computational issues arise when the problem has many design variables and many active constraints.

A variety of methods have been proposed to overcome these computational difficulties. Complex models typically call for more expensive analyses, leading to the use of approximation methods. A popular choice is to use an algorithm that generates a response surface for the problem, reducing the time for structural analysis [2,3].

Aerospace vehicle design is commonly constrained by buckling issues. The buckling of thin-walled shell structures often governs the design of light-weight aerospace vehicles. This applies to both spacecraft during launch and stage separation, as well as airplanes. A typical shell segment contains internal axial stiffeners and circumferential stiffeners – stringers and rings, respectively. Including buckling constraints in these problems make the optimization process more involved. At each iteration the solver must solve both the pre-buckling analysis and eigenvalue problem [3-5].

This paper presents an investigation into the issues with standard optimization engines and the problem of optimizing a multi-component system under buckling constraints. The convergence issues lead to the development of an FEA-based structural optimization algorithm to overcome the problems associated with buckling constraints. Included in the study is an analysis of common optimization algorithms and how their methods may be beneficial or detrimental to converging on an optimum solution.

Buckling analysis details follow, and the subsequent sections detail our approach to solving this optimization problem. Section 2 discusses the technical approach, including the model and investigated algorithms. The resulting developed algorithm is presented in Section 3, followed with discussion in Section 4 and conclusions in Section 5.

**Buckling Analysis**

Buckling loads are critical loads where certain types of structures become unstable, leading to large displacements and potentially material failure. In modeling the stability of fuselage structures, the problem setup is a thin-walled cylindrical shell structure under axial compression. This loading condition applies to spacecraft fuselage during the launch, flight, and stage separation phases.

The typical behavior of a thick metal cylindrical shell in axial compression is axisymmetric deformation once the maximum load-carrying capacity is reached. This continues past this limit load until the bifurcation of equilibrium states initiates. Beyond this bifurcation point the shell will deform axisymmetric, and then deform both axisymmetric and non-axisymmetric – i.e. the total deformation state is nonlinear. The behavior of a thin-walled cylindrical shell is shown in Figure 1, where the bifurcation point occurs before the limit point is reached. In other words, non-axisymmetric buckles start to appear before axisymmetric collapse [6].

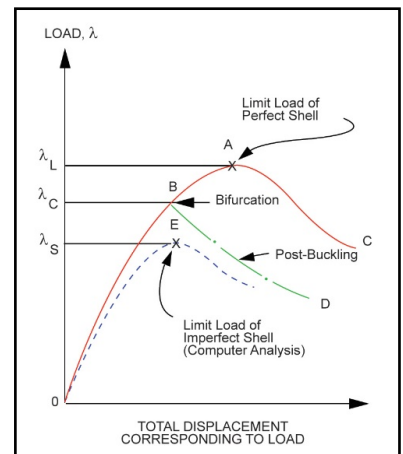

**Figure 1.** For a thin-walled cylinder in axial compression, bifurcation point, B, occurs before the limit point, A, is reached.

The bifurcation buckling in the *pre-buckling* phase has significant implications for FEA solvers. Another takeaway from Figure 1 is the maximum load-carrying capacity of the imperfect shell (dashed line) is below the load-carrying capacity of the perfect shell.

The FEA approach to solving buckling problems is two-fold: need to solve both static analysis and the eigenvalue problem at each optimization iteration. When bifurcation occurs after reaching the limit load, the pre-buckling solution is approximately independent of changes to member sections in the optimization process. The problem is then a relatively efficient gradient-based optimization because it's simple to calculate sensitivity of buckling loads.

However, for a multi-component shell or panel structure varying in thickness distribution the optimization is more complex. The distribution of pre-buckling stress becomes non-uniform, as in Figure 1. For a fuselage model, pre-buckling bending causes redistribution of stresses between the panel skin and the various segments of the stringers and rings. With multiple components there are two types of interaction between local & general buckling of the structures: (i) the general buckling mode for the system contains local components, and (ii) local buckling of the shell skin of a structure stiffened by ribs and stringers reduces the effective stiffness of that skin as it acts in a general buckling mode [3].

Eigenvalue buckling analysis predicts the theoretical critical buckling load – the bifurcation point – of an ideal linear elastic structure. Under a reference load $P_{ref}$, analysis of the stresses gives the stress-stiffness matrix $K_{ss}$. Buckling loads are then calculated by solving the eigenvalue problem:

$$([K] - \lambda_j [K_{ss}])\{u_j\} = \{0\}$$

where $[K]$ is the global stiffness matrix, $[K_{ss}]$ is the stress-stiffness (global geometric) matrix, and $\lambda_j$ and $\{u_j\}$ are the $j$th eigenvalue and eigenvector.

The first eigenvalue, $\lambda_1$, is the critical eigenvalue associated with buckling. The critical or buckling load is $P_{cr}$, associated with this eigenvalue. When a unit reference load is applied to the model, the limit load $P_{cr}$ equals the first mode eigenvalue $\lambda_1$.

## II. TECHNICAL APPROACH

The approach in this study was to apply standard optimization algorithms in order to provide insight to the problem – design variable trends, convergence issues, failure modes, etc. The hope was that these insights would offer a clear picture of the true issues, and provide helpful characteristics and techniques in designing a new optimization algorithm for the buckling problem.

### Model Analysis and Optimization

Design variables are system parameters that can be changed to improve the system performance. The model used in this study has five design variables, which are component sections of the fuselage panel. The modifiable parameters are component thicknesses. Figure 2 shows the fuselage panel model used in this study, with the five discretized sections color coded. The model geometry is representative of the technical problems in this study: multiple components for local optimizations, panel with stringers and rings, and varying thickness. It's also a useful model for the problem because an optimized solution isn't achieved with the solvers, as discussed in the next section.

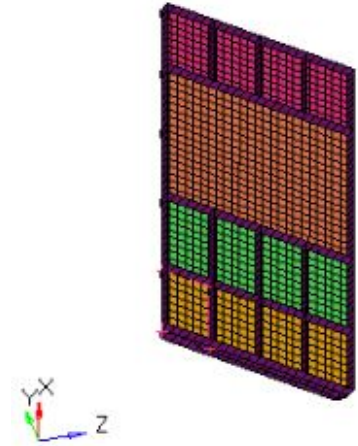

**Figure 2.** Fuselage panel model.

The panel is constrained in the vertical axis, and loaded with unit pressure across the top of the panel. This represents *quasi-static loading*. A static load is time independent, and a dynamic load is time dependent – for which inertial effects cannot be ignored. A quasi-static load is time dependent but is "slow enough" such that inertial effects can be ignored, and buckling analysis can be performed in this manner. Baseline analysis of the panel model is shown in the left side of Figure 3. The panel is to be optimized for minimal mass, where the component sections can vary in thickness. The optimization problem can be formally defined as

Minimize $F(x)$,
subject to
$G_i(x) = 0, i = 1, \dots, k_e$
$G_i(x) \leq 0, i = k_e + 1, \dots, k$
$lb \leq x \leq ub$

where the objective function $F$ is the system mass. In this problem $x$ is the vector of component thicknesses throughout the model. The constraints $G$ are mass greater than zero and $\lambda \geq 1.30$. The design space is defined by the bounds on component thicknesses, which are 0.001 to 1.0.

Standard gradient-based optimization (Nastran) of the model is shown in the right side of Figure 3. The important takeaway here is the difference in shape factors from the first buckling mode. The baseline analysis shows the full model approaching buckling in unison. That is, with a constant thickness across model components, there's a continuous buckling shape. The "optimized" model on the right of Figure 3, however, shows individual buckling shapes for the individual sections segmented by ribs. Thus, discretized sections of the model each have unique buckling responses.

This illustrates the effects of stiffness components (stringers and ribs) on the full model in buckling. First, the general buckling mode for the system contains local components. Second, local buckling of the shell skin of a structure stiffened by ribs and stringers reduces the effective stiffness of that skin as it acts in a general buckling mode [3].

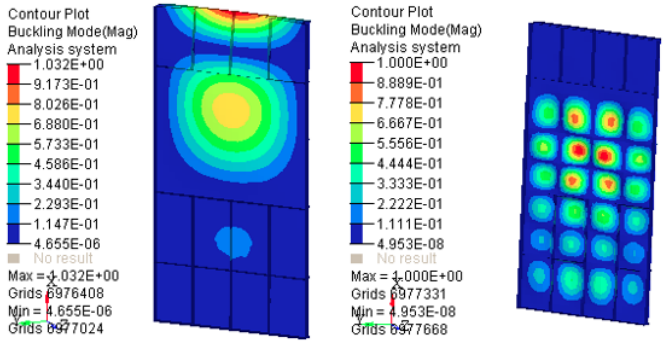

**Figure 3.** Buckling mode one shapes for the baseline analysis run (Left??) compared to an unconverged solution from the Nastran solver (right).

**Algorithms**

The Nastran solver is a gradient-based method, numerically searching for an optimum solution. At a given point in the design space, the gradients of the objective function and constraints are determined. This gives the solver a search direction because the gradient vector points in the direction of increasing the objective function. The process takes a number of steps in this search direction, which is equivalent to a number of function evaluations in numerical optimization, resulting in a new design. From this new point in design space, the process of computing gradients and establishing another search direction is repeated. This continues until the design converges on a solution – no further improvement in the objective function without violating any of the constraints [5].

The Nastran optimization run shown on the right side of Figure 3 was unable to converge on a solution. Convergence was achieved only when more iterations were added to the optimization. This is okay for simple models, but can become computationally burdensome for larger models with more design variables and/or constraints. Yet the unconverged solution is still valuable for this study, serving as a 'baseline optimum' for comparison (see Table 1).

| Mass | 4.123 |
|---|---|
| $\lambda_1$ | 1.293 |
| $T_1$ | 2.820E-02 |
| $T_2$ | 2.857E-02 |
| $T_3$ | 1.486E-02 |
| $T_4$ | 2.856E-02 |
| $T_5$ | 2.708E-02 |

**Table 1.** Nastran "baseline optimum" results.

The desire for computational efficiency in optimization of large system models has led to the development of approximation methods such as response surfaces and derivative-based local approximations. Global approximation methods (response surface based) are designed to be efficient and hence they are preferred methods when dealing with large design spaces, non-linear responses, or multi-objective optimization problems. Global optimization methods use higher order polynomials to approximate the original structural optimization problem over a wide range of design variables [3].

Response surface approximations help work around computational costs of expensive analysis, such as a system with many design variables and constraints. The approximation model then replaces the analysis program in the optimization. By using an approximation of the model, however, accuracy may be compromised for the higher efficiency. In the Adaptive Response Surface Method (ARSM) approach, the objective and constraint functions are approximated in terms of design variables using a second order polynomial [5]. In Altair HyperStudy, this is the default optimization engine for single-objective problems. Attempting to optimize the fuselage panel resulted in runs like that shown in Figure 4. The solver iteratively adjusts the components' thicknesses, then decreases them uniformly before 'getting stuck' at very small panel thicknesses. That is, the ARSM solver moves into the infeasible design space and is unable to recover.

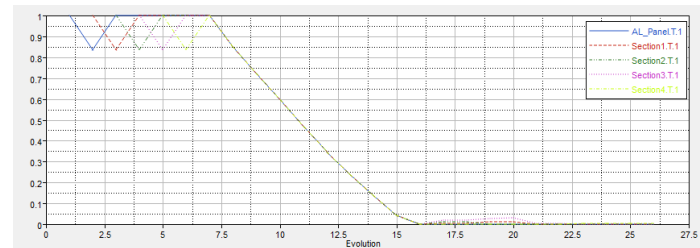

**Figure 4.** Typical optimization evolution for response surface approximation methods – ARSM and GRSM.

The Global Response Surface Method (GRSM) is a preferred response surface algorithm for large design spaces. Additional designs are generated globally, so the solver looks to optimize both locally and globally, in parallel [5]. Like the ARSM results, the approximation method falls into the infeasible design space and doesn't recover.

This behavior was also observed for runs with the Sequential Quadratic Programming (SQP) algorithm, despite being perhaps the optimization method of choice across industries. The fundamental principle behind SQP is to create a quadratic approximation of the Lagrangian and to solve that quadratic problem to define the search direction. Then like other solvers, the design variables are modified following that direction, improving the objective function without violating constraints. Theoretically SQP is more accurate but not as efficient as the response surface methods [5]. Thus one would expect the results, where SQP failed to come close to converging on a feasible design.

Convergence was achieved, however, with the Method of Feasible Directions (MFD) optimization engine. The fundamental principle behind the MFD is to progress from one feasible design to another, reducing the objective function. The key characteristic is the constraints at the new design point should never be violated. The method is decomposed into two sub-problems: (i) find a search direction such that the objective function reduces and the search remains in the feasible domain, and (ii) determine the step size [5]. Never dipping into the infeasible design space, MFD is inherently a slower algorithm than the previously discussed methods. And the subsequent small improvements on the objective function from iteration-to-iteration are subject to easily satisfy the convergence criteria: absolute and relative objective-based convergence. Thus, the converged solutions are not guaranteed to be optimal. Furthermore, for a model of many components, the converged solution will almost certainly be suboptimal. The results are shown in Table 2.

The relative success of MFD on the panel model brings to light the advantage of staying in the feasible design space. A possible reason is that with a model of many components for optimizations, it's viable for a solver to minimize each design variable as much as possible, before considering the combined effect on buckling response. MFD avoids this behavior, whereas approximation methods are particularly susceptible.

The use of a genetic algorithm (GA) for the problem, as expected, was far from converging on a solution. The solver was stopped at 1400 generations (approximately 48 hours). With more than a few design variables, GAs need many generations to find the most fit members of the design space; GAs rely on improved solutions over successive generations, so a large number of design iterations are typically needed, especially for large system models. GAs do however offer some benefits for optimizing over a large solution space: they provide non-local solutions and are well suited for both discrete and continuous optimization problems [3]. The first two attributes can be useful in our problem of optimizing many components within a system model.

Across optimization methods, the algorithms were largely unable to correlate design variable change with the eigenvalue response. Along with the observation made from Figure 3 that discretized sections produce individual buckling shapes, the following hypothesis is made: If there are at least as many responses as there are design variables, then the solver will find a correlation. Even though the first buckling mode is the critical mode for stability, additional eigenvalue responses for the higher modes may help solvers recognize the relationship between design variable change and buckling response.

Figure 5 shows typical behavior for the GRSM constrained for buckling responses $\lambda_1 - \lambda_{25}$. Although convergence is still not achieved, the solver does not get 'stuck' in the infeasible design space, and a correlation for at least one design variable is starting to be recognized.

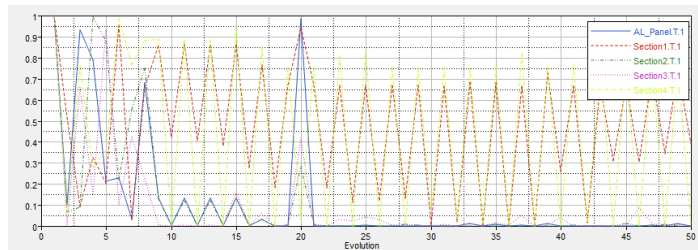

**Figure 5.** Typical plot of design variable vs. iteration for the GRSM optimizations under 25 buckling mode responses.

The increased number of buckling responses also improved the MFD converged solutions, as shown in Table 2. The trend is clear: more eigenvalue responses result in a better converged solution. Yet, as discussed earlier, these solutions are clearly suboptimal; they're inconsistent, and the Nastran baseline optimum is a better solution.

|  | **Objective** | **Response** |
|---|---|---|
| **1λ** | m=65.67kg | λcrit=3.385 |
| **5λ** | m=24.33kg | λcrit=3.127 |
| **15λ** | m=9.102kg | λcrit=6.027 |
| **25λ** | m=9.102kg | λcrit=6.027 |

**Table 2.** Converged results for the Method of Feasible Directions for a varied number of eigenvalue response modes.

Two consistent themes were observed for the tested optimization algorithms. First, they all struggled to correlate changes in component thickness to the eigenvalue responses. Specifically, approximation methods struggle to handle many design variables under buckling constraints. The second observation is this issue was slightly mitigated with more eigenvalue responses. These lessons influence the design of a new optimization algorithm, *eigenOpt*.

## III. RESULTS

### eigenOpt

The relative success of MFD in optimizing the fuselage panel implies there may be useful characteristics of the algorithm. MFD stays away from the infeasible design space, although this feature may ultimately penalize computational efficiency.

Gradient-based algorithms, like Nastran, are local approximation methods. They are effective when the derivatives of the system responses with respect to design variables can be computed easily and inexpensively. This can present issues when the system responses are dependent on multiple changing design variables – multiple component thicknesses influencing eigenvalue responses. But the methodology works well for finding local optimums.

Building off this logic, eigenOpt is an algorithm developed for this problem. The fundamental concept of eigenOpt is each design variable is treated with its own unique response, within the scope of optimizing the full system model. That is, for a given component of the panel model, the thickness is modified (or held constant) depending on that component's eigenvalue response, relative to the other eigenvalue responses throughout the panel. The sensitivity of a component to buckling is the driving metric, which we call *buckling stability*.

Buckling stabilities $\beta_i^+$ and $\beta_i^-$ are calculated from the critical value of the load factor $\lambda_i$ for each section $i$, quantifying the susceptibility of a given section to buckling. The derivation follows:

$$([K] - \lambda_j[K_{ss}])\{u_j\} = \{0\}$$

where $[K]$ is the global stiffness matrix, $[K_{ss}]$ is the stress-stiffness (global geometric) matrix, and $\lambda_j$ and $\{u_j\}$ are the $j$th eigenvalue and eigenvector. Each section has its own stiffness matrix $k_i$, which is a function of the section's thickness $t$. The global stiffness matrix is then composed of $k_{i-N}$ matrices for N total sections.

Following from earlier, the smallest positive eigenvalue $\lambda_1$ is the critical value of the load factor, or the smallest load at which buckling occurs. Solving for the eigenvalues gives

$$\lambda_j = \frac{\{u_j\}^T[K]\{u_j\}}{\{u_j\}^T[K_{ss}]\{u_j\}}$$

The quantity of interest for buckling stabilities is the change in eigenvalue response for a given change in thickness. Thus we differentiate the previous expression with respect to the design variable $t$:

$$\frac{\partial \lambda_j}{\partial t} = \frac{\{u_j\}^T[\frac{\partial[K]}{\partial t} - \frac{\lambda_j \partial[K_{ss}]}{\partial t}]\{u_j\}}{\{u_j\}^T[K_{ss}]\{u_j\}}$$

The denominator is constant for a change in thickness, so we ignore it and result with

$$\Delta \lambda_j = \{u_j\}^T([\Delta K] - \lambda_j[\Delta K_{ss}])\{u_j\}$$

The above equation is used with small changes in the thickness of a model section $i$ to give the buckling stability of that section. The change in $[\Delta K]$ is equal to the change in the section's stiffness matrix $[\Delta k_i]$.

The change in the global stress-stiffness matrix depends on the pre-buckling stress distributions, solved in the static analysis part of optimization. But with small changes in the thickness at each optimization iteration this is negligible. Thus, we ignore it in calculating the $j$th eigenvalue of the $i$th section:

$$\Delta \lambda_{ij} = \{u_j\}^T[\Delta k_i]\{u_j\}$$

It follows that when the thickness of a section $i$ is decreased,

$$[\Delta k_i] = [\Delta k_i]^- = [k_i(t - \Delta t)] - [k_i(t)]$$

And similarly for an increase in the section thickness,

$$[\Delta k_i] = [\Delta k_i]^+ = [k_i(t + \Delta t)] - [k_i(t)]$$

Finally we arrive at the buckling stability values, one each for thickness increase and decrease,

$$\beta_i^- = \{u_1\}^T[\Delta k_i]^-\{u_1\}$$
$$\beta_i^+ = \{u_1\}^T[\Delta k_i]^+\{u_1\}$$

The subscript of the eigenvector is 1 in these equations, corresponding with the smallest load at which buckling occurs. Thus the objective is to increase the thickness of sections with high values of $\beta_i^+$ and decrease the thickness of sections with high values of $\beta_i^-$. The optimization steps are

1. Problem setup: mesh the primary structure, discretize model into $i$ sections of different thicknesses, setup variables and constraints, set objective to minimize mass.
2. Run stress analysis, solving for buckling eigenvalues
3. Calculate buckling stabilities $\beta_i^+$ and $\beta_i^-$ (see below) for each discretized section.
4. For sections with highest values of $\beta_i^{+/-}$, change the thickness: increase for high values of $\beta_i^+$ and decrease for high values of $\beta_i^-$. Modifications are subject to the size constraint for each section where $t_{section\,min} \leq t_{section} \leq t_{section\,max}$.
5. Repeat steps 2-4 until solution converges.

A problem-specific threshold needs to be set for the number of changing model sections in each iteration. For example, setting this threshold at 0.20 would change the thicknesses of the components with the 20% highest values of buckling sensitivity, positive or negative. It may be beneficial to decrease this threshold with increasing iterations to help promote convergence.

The convergence criteria are to terminate at a maximum number of iterations, and relative objective convergence. An objective-dependent criterion is chosen because otherwise one can imagine the likely scenario where the solution is close to an optimum, but the thickness of two (or more) sections fluctuate back and forth without changing the mass. A design variable-dependent criterion would carry on forever.

## IV. DISCUSSION & CONCLUSIONS

Some fundamental characteristics of eigenOpt are pulled from the investigated algorithms. An MFD-like feature is that eigenOpt aims to stay in the feasible design space. By only modifying a select portion of system components each iteration, and also increasing thicknesses when necessary, the optimization is guarded against the infeasible design space.

eigenOpt is a *pseudo-evolutionary algorithm* for several reasons. First, the method uses both exploitation and exploration, looking for solutions locally each iteration, but also globally overall. It does this by calculating the buckling stabilities across the system, and modifying locally on the more significant components. Second, the algorithm evolves to a lower "mutation rate", i.e. close to convergence the relative modification in design variables decreases. This again plays into the focus on staying in the feasible design space. And another feature making eigenOpt an evolutionary-like algorithm is the "fitness ranking", where the change in design variables each iteration is based on the "most fit" members, i.e. the components with greatest buckling sensitivities.

Like evolutionary algorithms, this method may not ensure a theoretical optimum, but is designed to efficiently produce a reasonable discrete solution.

The next step in this project is to test the eigenOpt algorithm, but in an FEA package other than HyperWorks. Although HyperWorks has functionality to incorporate an external optimization engine, the deck is not accessible mid-optimization. Thus, the requisite component-specific stiffness matrices, or eigenvalues for that matter, cannot be extracted for use in eigenOpt.

Future work would consider the influence of imperfections on the model. The load-carrying capacity of optimally designed structures is often reduced more by imperfections than is so for other structures [7]. An imperfect shell structure will bend when any loads are applied. This pre-buckling bending causes redistribution of stresses between the panel skin and the various segments of the stringers and rings [6]. Fuselage structures may buckle due to launch loads, but also combined with circumferentially varying dynamic pressure. In future work this could be an additional optimization constraint.